\documentstyle[12pt]{article}
\newtheorem{defn}{{\sc Definition}}[section]
\newtheorem{thm}{Theorem}[section]
\newtheorem{Theorem}{Theorem}[section]
\newtheorem{lem}[thm]{Lemma}

\newtheorem{cor}[thm]{Corollary}
\newcommand{\sthm}{\begin{Theorem}}         
\newcommand{\ethm}{\end{Theorem}}           
\newtheorem{Corollary}[Theorem]{Corollary}   
\newcommand{\scor}{\begin{Corollary}}       
\newcommand{\ecor}{\end{Corollary}}         
\newcommand{\pf}{ \par \vspace{1ex} \noindent {\sc Proof.} \hspace{2mm}}

\begin{document}

\title{Steady State Coexistence Solutions of Reaction-Diffusion Competition Models}
\author{Joon Hyuk Kang\\
Department of Mathematics\\
Andrews University\\
Berrien Springs, MI. 49104\\
U.S.A.\\
kang@andrews.edu}
\date{ }
\maketitle
\vskip30pc
\begin{abstract}
Two species of animals are competing in the same environment. Under what conditions do they coexist peacefully? Or under what conditions does either one of the two species become extinct, that is, is either one of the two species excluded by the other? It is natural to say that they can coexist peacefully if their rates of reproduction and self-limitation are relatively larger than those of competition rates. In other words, they can survive if they interact strongly among themselves and weakly with others. We investigate this phenomena in mathematical point of view.
\\[0.1in] 
In this paper, we concentrate on coexistence solutions of the competition model 
$$\left\{ \begin{array}{l} 
\left.\begin{array}{l}
\Delta u + u(a - g(u,v)) = 0\\ 
\Delta v + v(d - h(u,v)) = 0
\end{array} \right.\;\;\mbox{in}\;\;\Omega,\\
u|_{\partial\Omega} = v|_{\partial\Omega} = 0,
\end{array} \right.$$
This system is the general model for the steady state of a competitive interacting system. The techniques used in this paper are elliptic theory, super-sub solutions, maximum principles, implicit function theorem and spectrum estimates. The arguments also rely on some detailed properties for the solution of logistic equations. 
\end{abstract}

\section{Introduction}
A lot of research has been focused on reaction-diffusion equations modeling of various systems in mathematical biology, especially the elliptic steady states of competitive and predator-prey interacting processes with various boundary conditions. In the earlier literature, investigations into mathematical biology models were concerned with studying those with homogeneous Neumann boundary conditions. From here on, the more important Dirichlet problems, which allow flux across the boundary, became the subject of study.
\\[0.1in]
Suppose two species of animals, rabbits and squirrels for instance, are competing in a bounded domain $\Omega$. Let $u(x,t)$ and $v(x,t)$ be densities of the two habitats in the place $x$ of $\Omega$ at time $t$. Then we have the following biological interpretation of terms.
\\[0.1in]
$(A)$ The partial derivatives $u_{t}(x,t)$ and $v_{t}(x,t)$ mean the rate of change of densities with respect to time $t$.\\
$(B)$ The laplacians $\Delta u(x,t)$ and $\Delta v(x,t)$ imply the diffusion or migration rates.\\
$(C)$ The rates of self-reproduction of each species of animals are expressed as multiples of some positive constants $a, d$ and current densities $u(x,t), v(x,t)$, i.e. $au(x,t)$ and $dv(x,t)$ which will increase the rate of change of densities in $(A)$, where $a > 0, d > 0$ are called self-reproduction constants.\\
$(D)$ The rates of self-limitation of each species of animals are multiples of some positive constants $b, f$ and the frequency of encounters among themselves $u^{2}(x,t), v^{2}(x,t)$, i.e. $bu^{2}(x,t)$ and $v^{2}(x,t)$ which will decrease the rate of change of densities in $(A)$, where $b > 0, f > 0$ are called self-limitation constants.\\
$(E)$ The rates of competition of each species of animals are multiples of some positive constants $c, e$ and the frequency of encounters of each species with the other $u(x,t)v(x,t)$, i.e. $cu(x,t)v(x,t)$ and $eu(x,t)v(x,t)$ which will decrease the rate of change of densities in $(A)$, where $c > 0, e > 0$ are called competition constants.\\
$(F)$ We assume that both species of animals are not staying on the boundary of $\Omega$.
\\[0.1in]
Combining all those together, we have the following dynamic model
$$\left\{ \begin{array}{l} 
\left.\begin{array}{l}
u_{t}(x,t) = \Delta u(x,t) + au(x,t) - bu^{2}(x,t) - cu(x,t)v(x,t)\\ 
v_{t}(x,t) = \Delta v(x,t) + dv(x,t) - fv^{2}(x,t) - eu(x,t)v(x,t)
\end{array} \right.\;\;\mbox{in}\;\;\Omega \times [0,\infty),\\
u(x,t) = v(x,t) = 0\;\;\mbox{for}\;\;x \in \partial\Omega,
\end{array} \right.$$
or equivalently,
$$\left\{ \begin{array}{l} 
\left.\begin{array}{l}
u_{t}(x,t) = \Delta u(x,t) + u(x,t)(a - bu(x,t) - cv(x,t))\\ 
v_{t}(x,t) = \Delta v(x,t) + v(x,t)(d - fv(x,t) - eu(x,t))
\end{array} \right.\;\;\mbox{in}\;\;\Omega \times [0,\infty),\\
u(x,t) = v(x,t) = 0\;\;\mbox{for}\;\;x \in \partial\Omega,
\end{array} \right.$$
Here, we are interested in the time independent, positive solutions, i.e. the positive solutions $u(x), v(x)$ of
$$\left\{ \begin{array}{l} 
\left.\begin{array}{l}
\Delta u(x) + u(x)(a - bu(x) - cv(x)) = 0\\ 
\Delta v(x) + v(x)(d - fv(x) - eu(x)) = 0
\end{array} \right.\;\;\mbox{in}\;\;\Omega,\\
u|_{\partial\Omega} = v|_{\partial\Omega} = 0,
\end{array} \right.$$
which are called coexistence state or steady state. The coexistence state is the positive density solution depending only on spatial variable $x$, not on time variable $t$, and so the existence of that means the two species of animals can live peacefully and forever.
\\[0.1in]
A lot of work about the existence and uniqueness of coexistence state of the above steady state model has already been established during the last decade.(see ~\cite{cc87}, ~\cite{cc89}, ~\cite{cl84}, ~\cite{gl94}, ~\cite{gp91}, ~\cite{ll91}, ~\cite{l80}.) 
\\[0.1in]
In this paper, we study rather general types of system. We concern the existence and uniqueness of positive coexistence when the relative growth rates are nonlinear, more precisely, the existence and uniqueness of positive steady state of
$$\left\{ \begin{array}{l} 
\left.\begin{array}{l}
\Delta u + u(a - g(u,v)) = 0\\ 
\Delta v + v(d - h(u,v)) = 0
\end{array} \right.\;\;\mbox{in}\;\;\Omega,\\
u|_{\partial\Omega} = v|_{\partial\Omega} = 0,
\end{array} \right.$$
where $a, d$ are positive constants, $g, h$ are $C^{1}$ functions, $\Omega$ is a bounded domain 
in $R^{n}$ and $u, v$ are densities of two competitive species. 
\\[0.2in]
The followings are questions raised in the general model with nonlinear growth rates.
\\[0.1in]
Problem 1 : Under what conditions do they coexist? Under what conditions do they have unique steady state? When does either one of the species become extinct?
\\[0.2in]
Problem 2 : Assuming that they can coexist and the coexistence state is unique at a fixed self-reproduction $(a,d)$, can they still coexist regardless of slight change of that self-reproduction? 
\\[0.2in]
Problem 3 : This is a generalization of Problem 2. If we have the existence and uniqueness of coexistence state on the left boundary of a closed convex region $\Gamma$ for the reproduction $(a,d)$, can we extend the region $\Gamma$ to an open set including $\Gamma$ without losing the uniqueness?  
\\[0.2in]
\\[0.1in]
In section 3, some sufficient conditions to guarantee the existence and uniqueness of positive solutions were obtained, and we could also see that there is no positive solution for small self-reproduction rates mainly using upper-lower solutions and spectrum estimates, which solves Problem 1. In sections 4 and 5, we answer problems 2 and 3 using elliptic theory, maximum principles and implicit function theorem.       

\section{Preliminaries}
In this section we will state some mathematical preliminary results which will be useful for our later arguments.
\begin{defn}\label{defn:2.1}(Super and sub solutions)
The vector functions $(\bar{u}^{1},...,\bar{u}^{N})$ and $(\underline{u}^{1},...,\underline{u}^{N})$ form an super/sub solution pair for the system
$$\left\{\begin{array}{r}
\Delta u^{i} + g^{i}(u^{1},...,u^{N}) = 0\;\;\mbox{in}\;\;\Omega\\
u^{i} = 0\;\;\mbox{on}\;\;\partial\Omega
\end{array}\right.$$
if for $i = 1,...,N$
$$\left\{\begin{array}{r}
\Delta \bar{u}^{i} + g^{i}(u^{1},...,u^{i-1},\bar{u}^{i},u^{i+1},...,u^{N}) \leq 0\\
\Delta \underline{u}^{i} + g^{i}(u^{1},...,u^{i-1},\underline{u}^{i},u^{i+1},...,u^{N}) \geq 0\\
\mbox{in}\;\;\Omega\;\;\mbox{for}\;\;\underline{u}^{j} \leq u^{j} \leq \bar{u}^{j}, j \neq i,
\end{array}\right.$$
and
$$\left.\begin{array}{l}
\underline{u}^{i} \leq \bar{u}^{i}\;\;\mbox{on}\;\;\Omega\\
\underline{u}^{i} \leq 0 \leq \bar{u}^{i}\;\;\mbox{on}\;\;\partial\Omega.
\end{array}\right.$$
\end{defn}

\begin{lem}\label{lem:7.1}
If $g^{i}$ in the Definition \ref{defn:2.1} are in $C^{1}$ and the system admits an super/sub solution pair $(\underline{u}^{1},...,\underline{u}^{N}), (\bar{u}^{1},...,\bar{u}^{N})$, then there is a solution of the system in \ref{defn:2.1} with $\underline{u}^{i} \leq u^{i} \leq \bar{u}^{i}$ in $\bar{\Omega}$. If 
$$\left.\begin{array}{l}
\Delta\bar{u}^{i} + g^{i}(\bar{u}^{1},...,\bar{u}^{N}) \neq 0,\\
\Delta\underline{u}^{i} + g^{i}(\underline{u}^{1},...,\underline{u}^{N}) \neq 0
\end{array}\right.$$
in $\Omega$ for $i = 1,...,N$, then $\underline{u}^{i} < u^{i} < \bar{u}^{i}$ in $\Omega$.
\end{lem}

\begin{lem}\label{lem:7.2}(The first eigenvalue)
\begin{equation}
\left\{ \begin{array}{l}
-\Delta u + q(x)u = \lambda u\;\;\mbox{in}\;\;\Omega,\\
u|_{\partial \Omega} = 0,
\end{array} \right. \label{eq:pb}
\end{equation}
where $q(x)$ is a smooth function from $\Omega$ to $R$ and $\Omega$ is a bounded domain in $R^{n}$.\\
$(A)$ The first eigenvalue $\lambda_{1}(q)$ of (\ref{eq:pb}), denoted by simply $\lambda_{1}$ when $q \equiv 0$, is simple with a positive eigenfunction.\\
$(B)$ If $q_{1}(x) < q_{2}(x)$ for all $x \in \Omega$, then $\lambda_{1}(q_{1}) < \lambda_{1}(q_{2})$.\\
$(C)$(Variational Characterization of the first eigenvalue)\\
$$\lambda_{1}(q) = \min_{\phi \in W_{0}^{1}(\Omega),\phi \neq 
0}\frac{\int_{\Omega}(|\nabla \phi|^{2}+q\phi^{2})dx}{\int_{\Omega}\phi^{2}dx}$$.
\end{lem}

\begin{lem}\label{lem:7.3}
$$Lu = \sum_{i,j=1}^{n}a_{ij}(x)D_{ij}u + \sum_{i=1}^{n}a_{i}(x)D_{i}u + a(x)u = f(x)\;\;\mbox{in}\;\;\Omega,$$
where $\Omega$ is a bounded domain in $R^{n}$ and\\
(M1) $\partial \Omega \in C^{2,\alpha} (0 < \alpha < 1)$,\\
(M2) $|a_{ij}(x)|_{\alpha}, |a_{i}(x)|_{\alpha}, |a(x)|_{\alpha} \leq M (i,j = 1,...,n)$,\\
(M3) $L$ is uniformly elliptic in $\bar{\Omega}$, with ellipticity constant $\gamma$, i.e., for every $x \in \bar{\Omega}$ and every real vector $\xi = (\xi_{1},...,\xi_{n})$
$$\sum_{i,j=1}^{n}a_{ij}(x)\xi_{i}\xi_{j} \geq \gamma\sum_{i=1}^{n}|\xi_{i}|^{2}.$$
$1.$ Maximum principles\\
Let $u \in C^{2}(\Omega) \cap C(\bar{\Omega})$ be a solution of $Lu \geq 0 (Lu \leq 0)$ in $\Omega$.\\
$(A)$ If $a(x) \equiv 0$, then $\max_{\bar{\Omega}}u = \max_{\partial\Omega}u (\min_{\bar{\Omega}}u = \min_{\partial\Omega}u)$.\\
$(B)$ If $a(x) \leq 0$, then $\max_{\bar{\Omega}}u \leq \max_{\partial\Omega}u^{+} (\min_{\bar{\Omega}}u \geq -\max_{\partial\Omega}u^{-})$,\\
where $u^{+} = \max(u,0), u^{-} = -\min(u,0)$.\\
$(C)$ If $a(x) \equiv 0$ and $u$ attains its maximum (minimum) at an interior point of $\Omega$, then $u$ is identically a constant in $\Omega$.\\
$(D)$ If $a(x) \leq 0$ and $u$ attains a nonnegative maximum (nonpositive minimum) at an interior point of $\Omega$, then $u$ is identically a constant in $\Omega$.
\\[0.1in]
$2.$ Schauder's estimate\\
If $u \in C^{2,\alpha}(\bar{\Omega})$ and $u|_{\partial\Omega} = \phi \in C^{2,\alpha}(\partial\Omega)$, then
$$|u|_{2,\alpha} \leq c(|Lu|_{\alpha} + |u|_{0} + |\phi|_{2,\alpha}^{\partial\Omega}),$$
where the constant $c > 0$ is independent of $u$.
\end{lem}

\begin{lem}\label{lem:7.4}(Implicit Function Theorem)\\
Let $X, Y, Z$ be Banach spaces. For a given $(u_{0},v_{0}) \in X \times Y$ and $a, b > 0$, let $S = \{(u,v) : \parallel u - u_{0} \parallel \leq a, \parallel v - v_{0} \parallel \leq b\}.$ Suppose $F : S \rightarrow Z$ satisfies the following:\\
$(A)$ $F$ is continuous.\\
$(B)$ $F_{v}(\cdot,\cdot)$ exists and is continuous in $S$(in the operator norm).\\
$(C)$ $F(u_{0},v_{0}) = 0$.\\
$(D)$ $[F_{v}(u_{0},v_{0})]^{-1}$ exists and is a continuous map from $Z$ to $Y$.\\
Then there are neighborhoods $U$ of $u_{0}$ and $V$ of $v_{0}$ such that the equation $F(u,v) = 0$ has exactly one solution $v \in V$ for every $u \in U$ and the solution $v$ depends continuously on $u$.
\end{lem}  

We also need some information on the solutions of the following logistic equations.
\begin{lem}\label{lem:7.6}(in \cite{ll91})
$$\left\{ \begin{array}{l}
\Delta u + uf(u) = 0\;\; \mbox{in}\;\; \Omega,\\
u|_{\partial\Omega} = 0, u > 0,
\end{array} \right.$$
where $\Omega$ is a bounded domain in $R^{n}$ and\\ 
$(A)$ $f$ is a strictly decreasing $C^{1}$ function,\\
$(B)$ there exists $c_{0} > 0$ such that $f(u) \leq 0$ for $u \geq c_{0}$.
\\[0.1in]
$(1)$ If $f(0) > \lambda_{1}$, where $\lambda_{1}$ is the first eigenvalue of $-\Delta$ with homogeneous boundary condition, then the above equation has a unique positive solution. 
$(2)$ If $f(0) < \lambda_{1}$, then $u \equiv 0$ is the only nonnegative solution of the above equation.
\end{lem}

In the case $(1)$, we denote this unique positive solution as $\theta_{f}$.
The main property about this positive solution is that $\theta_{f}$ is larger as $f$ is larger, i.e. $\theta_{g} \leq \theta_{f}$ if $g \leq f$.(Proposition \ref{pro:6.1})

\section{Existence, Nonexistence and Uniqueness of steady state}
We consider the elliptic system
\begin{equation}
\left\{ \begin{array}{l} 
\left. \begin{array}{l}
\Delta u + u(a - g(u,v)) = 0\\
\Delta v + v(d - h(u,v)) = 0
\end{array} \right. \mbox{in $\Omega$,}\\
u|_{\partial \Omega} = v|_{\partial \Omega} = 0.
\end{array} \right. \label{eq:ia}
\end{equation}  
Here $\Omega$ is a bounded, smooth domain in $R^{n}$ and\\
$(U1)$ $g, h \in C^{1}$ are strictly increasing functions with respect to $u, v$,\\
$(U2)$ there exist $k_{1}, k_{2} > 0$ such that $g(u,0) > a$ for $u \geq k_{1}$ and\\
$h(0,v) > d$ for $v \geq k_{2}$.  
\\[0.1in]
If there is no competition between the species, that is, if we consider
$$\left\{\begin{array}{l}
\left.\begin{array}{rrl}
\Delta u + u(a - g(u,0)) & = & 0\\
\Delta v + v(d - h(0,v)) & = & 0
\end{array}\right.\;\;\mbox{in}\;\;\Omega,\\
u = v = 0\;\;\mbox{on}\;\;\partial\Omega,
\end{array}\right.$$
then by the Lemma \ref{lem:7.6}, the condition $a > \lambda_{1}, d > \lambda_{1}$(i.e. reproductions are relatively large.) were sufficient to guarantee the positive density solution $\theta_{a-g(\cdot,0)}, \theta_{d-h(0,\cdot)}$. But, if there is some competition between them, then as see in the following Theorem \ref{thm:2.1}, we should have larger lower bound for reproduction rates $a$ and $d$, i.e. we have stronger conditions $a > \lambda_{1} + g(0,k_{2}), d > \lambda_{1} + h(k_{1},0)$ to guarantee their coexistence.(i.e. the reproductions should be much larger.)
\\[0.1in]
The following is the main result:

\begin{thm}\label{thm:2.1} $(A)$ If $a > \lambda_{1} + g(0,k_{2})$ and $d > \lambda_{1} + 
h(k_{1},0)$, then (\ref{eq:ia}) has a positive solution $(u,v)$ with
$$\theta_{a-g(\cdot,k_{2})} < u < \theta_{a-g(\cdot,0)}, \theta_{d-h(k_{1},\cdot)} < v < \theta_{d-h(0,\cdot)}.$$
Conversely, any positive solution $(u,v)$ to (\ref{eq:ia}) must satisfy the inequalities.\\  
$(B)$ If $a > \lambda_{1} + g(0,k_{2})$ and $d > \lambda_{1} + h(k_{1},0)$ and
\begin{eqnarray*}
\;\;4\inf_{B}(\frac{\partial g}{\partial u})\inf_{B}(\frac{\partial h}{\partial v}) 
& > & \sup\frac{\theta_{a-g(\cdot,0)}}{\theta_{d-h(k_{1},\cdot)}}(\sup(\frac{\partial 
g}{\partial v}))^{2} + \sup\frac{\theta_{dh(0,\cdot)}}{\theta_{ag(\cdot,k_{2})}}(\sup(\frac{\partial 
h}{\partial u}))^{2}\\ 
& + & 2\sup(\frac{\partial g}{\partial v}) \sup(\frac{\partial h}{\partial u}), 
\end{eqnarray*}
where $B = [0,k_{1}] \times [0,k_{2}]$, then (\ref{eq:ia}) has a unique coexistence state.\\
$(C)$ If $a \leq \lambda_{1}$ or $d \leq \lambda_{1}$, then (\ref{eq:ia}) does not have any positive solution. 
\end{thm}
Biologically, we can interpret the conditions in Theorem \ref{thm:2.1} as follows. The constants $a, d$ and functions $g, h$ describe how species 1 ($u$) and 2 ($v$) interact among themselves and with each other. Hence, the both conditions in $(A)$ and $(B)$ imply that species 1 interacts strongly among themselves and weakly with species 2. Similarly for species 2, they interact more strongly among themselves than they do with species 1. The inequalities in the conclusion $(A)$ imply that the densities with competitions($u$ and $v$) are less than those without competition.($\theta_{a-g(\cdot,0)}$ and $\theta_{d-h(0,\cdot)}$) Furthermore, $(C)$ says that if one of the species has small reproduction, then it may be extinct, which means that the two species can not coexist.  
\pf 
$(A)$ Let $\bar{u} = \theta_{a-g(\cdot,0)}, \bar{v} = \theta_{d-h(0,\cdot)}$. Then since $g$ is increasing, we have
$$\left.\begin{array}{ll}
& \Delta \bar{u} + \bar{u}(a - g(\bar{u},\bar{v}))\\
= & \Delta \bar{u} + \bar{u}(a - g(\bar{u},0) + g(\bar{u},0) - g(\bar{u},\bar{v}))\\
= & \bar{u}(g(\bar{u},0) - g(\bar{u},\bar{v})) < 0.
\end{array}\right.$$
Similarly, we have
$$\Delta \bar{v} + \bar{v}(d - h(\bar{u},\bar{v})) < 0.$$
So, $(\bar{u},\bar{v})$ is a super solution of (\ref{eq:ia}).\\
Let $\underline{u} = \theta_{a-g(\cdot,k_{2})}, \underline{v} = \theta_{d-h(k_{1},\cdot)}$.\\
Then by the Maximum Principle, we obtain
$$\left\{\begin{array}{l}
\underline{u} \leq \theta_{a-g(\cdot,0)} \leq k_{1},\\
\underline{v} \leq \theta_{d-h(0,\cdot)} \leq k_{2}.
\end{array}\right.$$
Since $g$ is increasing, we get
$$\left.\begin{array}{ll}
& \Delta \underline{u} + \underline{u}(a - g(\underline{u},\underline{v}))\\
= & \Delta \underline{u} + \underline{u}(a - g(\underline{u},k_{2}) + g(\underline{u},k_{2}) - g(\underline{u},\underline{v}))\\
= & \underline{u}(g(\underline{u},k_{2}) - g(\underline{u},\underline{v})) \geq 0.
\end{array}\right.$$
Similarly, we get
$$\Delta \underline{v} + \underline{v}(d - h(\underline{u},\underline{v})) \geq 0.$$
Therefore, $(\underline{u},\underline{v})$ is a lower solution of (\ref{eq:ia}).\\
Furthermore, $\underline{u} < \bar{u}$ and $\underline{v} < \bar{v}$ in $\Omega$ and $\underline{u} = \bar{u} = \underline{v} = \bar{v} = 0$ on $\partial\Omega$.\\
So, by Lemma \ref{lem:7.1}, (\ref{eq:ia}) has a solution $(u,v)$ with
$$\theta_{a-g(\cdot,k_{2})} < u < \theta_{a-g(\cdot,0)}, \theta_{d-h(k_{1},\cdot)} < v < \theta_{d-h(0,\cdot)}.$$
Suppose $(u,v)$ is a positive solution to (\ref{eq:ia}). 
By the Mean Value Theorem, there is $v^{*}$ such that\\
$$g(u,v) = g(u,0) + \frac{\partial g(u,v^{*})}{\partial v}v.$$
Then 
$$\Delta u + u(a - g(u,0)) = \frac{\partial g(u,v^{*})}{\partial v}uv > 0\;\;\mbox{in}\;\;\Omega.$$
Hence, $u$ is a subsolution to
$$\left\{\begin{array}{l}
\Delta z + z(a - g(z,0)) = 0\;\;\mbox{in}\;\;\Omega,\\
z|_{\partial\Omega} = 0.
\end{array}\right.$$
Any sufficiently large positive constant is a super solution to 
$$\left\{\begin{array}{l}
\Delta z + z(a - g(z,0)) = 0\;\;\mbox{in}\;\;\Omega,\\
z|_{\partial\Omega} = 0.
\end{array}\right.$$
Therefore, by the super-sub solution method, we have
\begin{equation}
\left. \begin{array}{l}
u \leq \theta_{a-g(\cdot,0)}.
\end{array} \right. \label{eq:ib}
\end{equation}
The same argument shows
\begin{equation}
\left. \begin{array}{l}
v \leq \theta_{d-h(0,\cdot)}.
\end{array} \right. \label{eq:ic}
\end{equation}
For sufficiently small $\epsilon > 0$,
$$\left.\begin{array}{ll}
& \epsilon\Delta\theta_{d-h(0,\cdot)} + \epsilon\theta_{d-h(0,\cdot)}(d - h(0,\epsilon\theta_{d-h(0,\cdot)}))\\
= & \epsilon[\Delta\theta_{d-h(0,\cdot)} + \theta_{d-h(0,\cdot)}(d - h(0,\epsilon\theta_{d-h(0,\cdot)}))]\\
> & \epsilon[\Delta\theta_{d-h(0,\cdot)} + \theta_{d-h(0,\cdot)}(d - h(0,\theta_{d-h(0,\cdot)}))]\\
= & 0\;\;\mbox{in}\;\;\Omega,
\end{array}\right.$$
and so $\epsilon\theta_{d-h(0,\cdot)}$ is a sub solution to
$$\left\{\begin{array}{l}
\Delta z + z(d - h(0,z)) = 0\;\;\mbox{in}\;\;\Omega,\\
z|_{\partial\Omega} = 0.
\end{array}\right.$$
Since $d - h(0,k_{2}) < 0, k_{2}$ is a super solution to 
$$\left\{\begin{array}{l}
\Delta z + z(d - h(0,z)) = 0\;\;\mbox{in}\;\;\Omega,\\
z|_{\partial\Omega} = 0.
\end{array}\right.$$
Hence, by the super-sub solution method again,
$$\theta_{d-h(0,\cdot)} \leq k_{2}.$$
So,
$$g(u,v) \leq g(u,\theta_{d-h(0,\cdot)}) \leq g(u,k_{2})$$
since $g(u,z)$ is increasing.
Therefore, 
$$\left.\begin{array}{lll}
\Delta u + u(a - g(u,k_{2})) & \leq & \Delta u + u(a - g(u,v))\\
& = & 0\;\;\mbox{in}\;\;\Omega.
\end{array}\right.$$
Hence, $u$ is a super solution to 
$$\left\{\begin{array}{l}
\Delta z + z(a - g(z,k_{2})) = 0\;\;\mbox{in}\;\;\Omega,\\
z|_{\partial\Omega} = 0.
\end{array}\right.$$
Let $\phi_{1}$ be the first eigenvector of
$$\left\{\begin{array}{l}
\Delta u + \lambda_{1}z = 0\;\;\mbox{in}\;\;\Omega,\\
z|_{\partial\Omega} = 0.
\end{array}\right.$$
Then for sufficiently small $\epsilon > 0$,
$$a - g(\epsilon\phi_{1},k_{2}) - \lambda_{1} > 0\;\;\mbox{in}\;\;\Omega,$$
and
$$\left.\begin{array}{lll}
\Delta(\epsilon\phi_{1}) + \epsilon\phi_{1}(a - g(\epsilon\phi_{1},k_{2})) & = & \epsilon[\Delta(\epsilon\phi_{1}) + \phi_{1}(a - g(\epsilon\phi_{1},k_{2}))]\\
& & \epsilon(\Delta\phi_{1} + \lambda_{1}\phi_{1}) = 0\;\;\mbox{in}\;\;\Omega.
\end{array}\right.$$
Consequently, $\epsilon\phi_{1}$ is a sub solution to 
$$\left\{\begin{array}{l}
\Delta z + z(a - g(z,k_{2})) = 0\;\;\mbox{in}\;\;\Omega,\\
z|_{\partial\Omega} = 0.
\end{array}\right.$$
Hence, by the super-sub solution method again,
\begin{equation}
\left. \begin{array}{l}
\theta_{a-g(\cdot,k_{2})} \leq u.
\end{array} \right. \label{eq:id}
\end{equation}
The same argument shows
\begin{equation}
\left. \begin{array}{l}
\theta_{d-h(k_{1},\cdot)} \leq v.
\end{array} \right. \label{eq:ie}
\end{equation}
From (\ref{eq:ib}) to (\ref{eq:ie}), we have
\begin{equation}
\theta_{a-g(\cdot,k_{2})} \leq u \leq \theta_{a-g(\cdot,0)},\;\;\\
\theta_{d-h(k_{1},\cdot)} \leq v \leq \theta_{d-h(0,\cdot)} 
\label{eq:if}
\end{equation}
Consequently, for any positive solution $(u,v)$ of (\ref{eq:ia}), the inequalities (\ref{eq:if}) hold.
\\[0.1in]
$(B)$ Suppose $(u_{1},v_{1})$ and $(u_{2},v_{2})$ are positive solutions to (\ref{eq:ia}).\\
Let $p = u_{1} - u_{2}$ and $q = v_{1} - v_{2}$. Then
\begin{eqnarray*}
\Delta p + (a - g(u_{1},v_{1}))p & = & \Delta u_{1} - \Delta u_{2} + 
(a - g(u_{1},v_{1}))(u_{1} - u_{2})\\
 & = & -\Delta u_{2} - (a - g(u_{1},v_{1}))u_{2}\\ 
& = & -\Delta u_{2} - u_{2}(a - g(u_{2},v_{2}) + g(u_{2},v_{2}) - g(u_{1},v_{1}))\\ 
& = & -u_{2}(g(u_{2},v_{2}) - g(u_{1},v_{1}))\\
& = & -u_{2}(g(u_{2},v_{2}) - g(u_{1},v_{2}) + g(u_{1},v_{2}) - g(u_{1},v_{1}))\\
& = & -u_{2}(\frac{\partial g(\tilde{x},v_{2})}{\partial u}(-p) + \frac{\partial g(u_{1},\bar{x})}{\partial v}(-q)\\
& = & u_{2}(p\frac{\partial g(\tilde{x},v_{2})}{\partial u} + q\frac{\partial g(u_{1},\bar{x})}{\partial v})\;\;\mbox{in}\;\;\Omega, 
\end{eqnarray*}
where $\tilde{x}, \bar{x}$ are from Mean Value Theorem depending on $u_{1}, u_{2}, v_{1}, v_{2}$. Hence,
\begin{equation}
\Delta p + (a - g(u_{1},v_{1}))p - u_{2}(p\frac{\partial g(\tilde{x},v_{2})}{\partial u} + q\frac{\partial g(u_{1},\bar{x})}{\partial v}) = 0\;\;\mbox{in}\;\;\Omega.\label{eq:ig}
\end{equation}
Similarly, we can get
\begin{equation}
\Delta q + (d - h(u_{2},v_{2}))q - v_{1}(p\frac{\partial h(\tilde{y},v_{1})}{\partial u} + q\frac{\partial h(u_{2},\bar{y})}{\partial v}) = 0\;\;\mbox{in}\;\;\Omega,\label{eq:ih}
\end{equation}
where $\tilde{y}, \bar{y}$ are from Mean Value Theorem depending on $u_{1}, u_{2}, v_{1}, v_{2}$.
Since $\lambda_{1}(a - g(u_{1},v_{1})) = 0$, by the Variational Characterization of the first eigenvalue, 
\begin{equation}
\int_{\Omega}z(-\Delta z - (a - g(u_{1},v_{1}))z)dx \geq 0 
\label{eq:ii}
\end{equation}
for any $z \in C^{2}(\bar{\Omega})$ and $z|_{\partial\Omega} = 0.$
The same argument shows that 
\begin{equation}
\int_{\Omega}w(-\Delta w - (d - h(u_{2},v_{2}))w)dx \geq 0 
\label{eq:ij}
\end{equation}
for any  $w \in C^{2}(\bar{\Omega})$ and 
$w|_{\partial\Omega} = 0.$
From (\ref{eq:ig}) and (\ref{eq:ih}), we have
$$\left\{ \begin{array}{l}
-p\Delta p - (a - g(u_{1},v_{1}))p^{2} + u_{2}p(p\frac{\partial g(\tilde{x},v_{2})}{\partial u} + q\frac{\partial g(u_{1},\bar{x})}{\partial v}) = 0\\  
-q\Delta q - (d - h(u_{2},v_{2}))q^{2} + v_{1}q(p\frac{\partial h(\tilde{y},v_{1})}{\partial u} + q\frac{\partial h(u_{2},\bar{y})}{\partial v}) = 0
\end{array} \right.\;\;\mbox{in $\Omega.$}$$
Using (\ref{eq:ii}) and (\ref{eq:ij}), we have
$$\int_{\Omega}[u_{2}p(p\frac{\partial g(\tilde{x},v_{2})}{\partial u} + q\frac{\partial g(u_{1},\bar{x})}{\partial v}) + v_{1}q(p\frac{\partial h(\tilde{y},v_{1})}{\partial u} + q\frac{\partial h(u_{2},\bar{y})}{\partial v})] \leq 0.$$
Hence,
$$\int_{\Omega}[u_{2}\frac{\partial g(\tilde{x},v_{2})}{\partial u}p^{2} + (u_{2}\frac{\partial g(u_{1},\bar{x})}{\partial v} + v_{1}\frac{\partial h(\tilde{y},v_{1})}{\partial u})pq + v_{1}\frac{\partial h(u_{2},\bar{y})}{\partial v}q^{2}] \leq 0.$$
Therefore, $p \equiv q \equiv 0$ if we can show that
$$(u_{2}\frac{\partial g(u_{1},\bar{x})}{\partial v} + v_{1}\frac{\partial h(\tilde{y},v_{1})}{\partial u})^{2} - 4u_{2}v_{1}\frac{\partial g(\tilde{x},v_{2})}{\partial u}\frac{\partial h(u_{2},\bar{y})}{\partial v} < 0\;\;\mbox{in}\;\;\Omega,$$
which is true if
$$\left.\begin{array}{l}
u_{2}^{2}(\frac{\partial g(u_{1},\bar{x})}{\partial v})^{2} + v_{1}^{2}(\frac{\partial h(\tilde{y},v_{1})}{\partial u})^{2} + 2u_{2}v_{1}\frac{\partial g(u_{1},\bar{x})}{\partial v}\frac{\partial h(\tilde{y},v_{1})}{\partial u}\\
- 4u_{2}v_{1}\frac{\partial g(\tilde{x},v_{2})}{\partial u}\frac{\partial h(u_{2},\bar{y})}{\partial v} < 0\;\;\mbox{in}\;\;\Omega.
\end{array}\right.$$
i.e.,
$$\left.\begin{array}{l}
4u_{2}v_{1}\frac{\partial g(\tilde{x},v_{2})}{\partial u}\frac{\partial h(u_{2},\bar{y})}{\partial v} > u_{2}^{2}(\frac{\partial g(u_{1},\bar{x})}{\partial v})^{2} + v_{1}^{2}(\frac{\partial h(\tilde{y},v_{1})}{\partial u})^{2}\\
+ 2u_{2}v_{1}\frac{\partial g(u_{1},\bar{x})}{\partial v}\frac{\partial h(\tilde{y},v_{1})}{\partial u}\;\;\mbox{in}\;\;\Omega,
\end{array}\right.$$
or
$$\left.\begin{array}{l}
4\frac{\partial g(\tilde{x},v_{2})}{\partial u}\frac{\partial h(u_{2},\bar{y})}{\partial v} > \frac{u_{2}}{v_{1}}(\frac{\partial g(u_{1},\bar{x})}{\partial v})^{2} + \frac{v_{1}}{u_{2}}(\frac{\partial h(\tilde{y},v_{1})}{\partial u})^{2}\\
+ 2\frac{\partial g(u_{1},\bar{x})}{\partial v}\frac{\partial h(\tilde{y},v_{1})}{\partial u}\;\;\mbox{in}\;\;\Omega,
\end{array}\right.$$
This is the case from the hypothesis in the theorem and (\ref{eq:if}), and so the uniqueness is proved.
\\[0.1in]
$(C)$ Assume $a \leq \lambda_{1}$. Suppose $(u,v)$ is a nonnegative solution to (\ref{eq:ia}). Then since $g$ is an increasing function with respect to $u$ and $v$, 
$$\left.\begin{array}{ll}
& \Delta u + u(a - g(u,0))\\
= & \Delta u + u(a - g(u,v) + g(u,v) - g(u,0))\\
= & u(g(u,v) - g(u,0)) \geq 0.
\end{array}\right.$$
Therefore, $u$ is a sub solution to
$$\left\{\begin{array}{l}
\Delta u + u(a - g(u,0)) = 0\;\;\mbox{in}\;\;\Omega,\\
u|_{\partial\Omega} = 0.
\end{array}\right.$$
Any constant larger than $k_{1}$ is a super solution to
$$\left\{\begin{array}{l}
\Delta u + u(a - g(u,0)) = 0\;\;\mbox{in}\;\;\Omega,\\
u|_{\partial\Omega} = 0.
\end{array}\right.$$
Hence, by the Lemma \ref{lem:7.1}, there is a solution $\bar{u}$ of
$$\left\{\begin{array}{l}
\Delta u + u(a - g(u,0)) = 0\;\;\mbox{in}\;\;\Omega,\\
u|_{\partial\Omega} = 0
\end{array}\right.$$
such that $0 \leq u \leq \bar{u}$. But, since $a \leq \lambda_{1}$, $\bar{u} \equiv 0$ by $(2)$ of Lemma \ref{lem:7.6}, and so $u \equiv 0$. 

\section{Uniqueness with small perturbation of reproduction rates}
We consider the model
\begin{equation}
\left\{ \begin{array}{l} 
\left. \begin{array}{l}
\Delta u +u(a - g(u,v))=0\\
\Delta v +v(d - h(u,v))=0
\end{array} \right. \mbox{in $\Omega$,}\\
u|_{\partial \Omega} = v|_{\partial \Omega} = 0.
\end{array} \right. \label{eq:ja}
\end{equation} 
Here $\Omega$ is a bounded, smooth domain in $R^{n}$ and\\
(P1) $g, h \in C^{1}$ are strictly increasing functions with respect to $u$ and $v$, and $g(0,0) = h(0,0) = 0$,\\
(P2) there are $k_{1}, k_{2} > 0$ such that $g(u,0) > a > \lambda_{1}$ for $u \geq k_{1}$ and $h(0,v) > d > \lambda_{1}$ for $v \geq k_{2}$.
\\[0.2in]
The following is the main theorem.
\begin{thm}\label{thm:3.1} Suppose\\
$(A)$ $a > \lambda_{1}(g(0,\theta_{d-h(0,\cdot)})), d > \lambda_{1}(h(\theta_{a-g(\cdot,0)},0))$,\\
$(B)$ $(\ref{eq:ja})$ has a unique coexistence state $(u,v)$,\\ 
$(C)$ the Frechet derivative of $(\ref{eq:ja})$ at $(u,v)$ is invertible.\\
Then there is a neighborhood $V$ of $(a,d)$ in $R^{2}$ such that if $(a_{0},d_{0}) \in V$, 
then $(\ref{eq:ja})$ with $(a,d) = (a_{0},d_{0})$ has a unique coexistence state.
\end{thm}
Theorem \ref{thm:3.1} looks like the consequence of Implicit Function Theorem. But the inverse function theorem only guaranteed the uniqueness locally. Theorem \ref{thm:3.1} concluded the global uniqueness. The techniques we will use includes naturally Implicit Function Theorem and a priori estimates on solutions of (\ref{eq:ja}).
\\[0.1in]
Biologically, the first condition in this theorem indicates that the rates of self-reproduction is large. The condition of invertibility of Frechet derivative also illustrates that the rates of self-limitation is relatively larger than those of competitions which will be in Theorem \ref{thm:3.2}. Then the conclusion says that small perturbation of reproduction rates does not lose the existence and uniqueness of positive steady state, i.e. they can still coexist peacefully even if there is some slight change of reproduction rates.
\pf Since the Frechet derivative of (\ref{eq:ja}) at $(u,v)$ is invertible, by the Implicit Function Theorem, there is a neighborhood $V$ of $(a,d)$ in $R^{2}$ and a neighborhood $W$ of $(u,v)$ in $[C_{0}^{2+\alpha}(\bar{\Omega})]^{2}$ such that for all $(a_{0},d_{0}) \in V$, there is a unique positive solution $(u_{0},v_{0}) \in W$ of (\ref{eq:ja}). Suppose the conclusion of the theorem is false. Then there are sequences $(a_{n},d_{n},u_{n},v_{n}), (a_{n},d_{n},u_{n}^{*},v_{n}^{*})$ in $V \times [C_{0}^{2+\alpha}(\bar{\Omega})]^{2}$ such that $(u_{n},v_{n})$ and $(u_{n}^{*},v_{n}^{*})$ are the positive solutions with  $(a,d) = (a_{n},d_{n})$ and $(u_{n},v_{n}) \neq (u_{n}^{*},v_{n}^{*})$ and $(a_{n},d_{n}) \rightarrow (a,d)$. By the standard elliptic theory, $(u_{n},v_{n}) \rightarrow (\bar{u},\bar{v})$ and   
$(u_{n}^{*},v_{n}^{*}) \rightarrow (u^{*},v^{*})$ in $C^{2,\alpha}$, and $(\bar{u},\bar{v})$, $(u^{*},v^{*})$ are solutions of (\ref{eq:ja}). Claim $\bar{u} > 0, \bar{v} > 0, u^{*} > 0, v^{*} > 0.$ It is enough to show that $\bar{u}$ and $\bar{v}$ are not identically zero because of the Maximum Principle. Suppose not, then by the Maximum Principle again, one of the 
following cases should hold: (1) $\bar{u}$ is identically zero and $\bar{v} > 0.$ (2) $\bar{u} > 0$ and $\bar{v}$ is identically zero. (3) $\bar{u}$ is identically zero and $\bar{v}$ is identically zero.\\
Without loss of generality, assume $\bar{u}$ is identically zero.\\
Let $\tilde{u_{n}} = \frac{u_{n}}{\parallel u_{n} \parallel_{\infty}},\tilde{v_{n}} = v_{n}$ for all $n\in N$. Then
$$\left\{ \begin{array}{l}
\Delta \tilde{u_{n}} + \tilde{u_{n}}(a_{n} - g(u_{n},\tilde{v_{n}})) = 0\\ 
\Delta \tilde{v_{n}} + \tilde{v_{n}}(d_{n} - h(u_{n},\tilde{v_{n}})) = 0 
\end{array} \right. \mbox{in $\Omega.$}$$
From the elliptic theory, $\tilde{u_{n}} \rightarrow \tilde{u}$ and
$$\left\{ \begin{array}{l}
\Delta \tilde{u} + \tilde{u}(a - g(0,\bar{v})) = 0\\
\Delta \bar{v} + \bar{v}(d - h(0,\bar{v})) = 0
\end{array} \right. \mbox{in $\Omega,$}$$
since $g, h$ are continuous. i.e., $a = \lambda_{1}(g(0,\bar{v}))$.\\ 
(1) If $\bar{v} \equiv 0$, then by the monotonicity of $g$ and $\lambda_{1}$, $a = \lambda_{1}(g(0,\bar{v})) = \lambda_{1}(g(0,0)) \leq \lambda_{1}(g(0,\theta_{d-h(0,\cdot)}))$ which contradicts our assumption.\\
(2) If $\bar{v}$ is not identically zero, then $\bar{v} = \theta_{d-h(0,\cdot)}$ and so $a = \lambda_{1}(g(0,\bar{v})) = \lambda_{1}(g(0,\theta_{d-h(0,\cdot)}))$ which is also a contradiction to our assumption.   
Consequently, $(\bar{u},\bar{v})$ and $(u^{*},v^{*})$ are coexistence states for $(a,d).$ But, since the coexistence state with respect to $(a,d)$ is unique, $(\bar{u},\bar{v})=(u^{*},v^{*})=(u,v)$. 
But, since $(u_{n},v_{n})\neq(u_{n}^{*},v_{n}^{*})$, it contradicts the Implicit Function Theorem.                           
\vskip1pc
The proof of the theorem also tells us that if one of the species becomes extinct, in other word, if one is excluded by others, then that means the reproduction rates are small, i.e. the region condition of reproduction rates $(A)$ is reasonable.
\begin{thm}\label{thm:8.1} If $(a_{n},d_{n},u_{n},v_{n}) \rightarrow (a,d,u,v)$ and if $u \equiv 0$ or $v \equiv 0$, then\\
$a \leq \lambda_{1}(g(0,\theta_{d-h(0,\cdot)}))$ or $d \leq \lambda_{1}(h(\theta_{a-g(\cdot,0)},0))$.
\end{thm}  
The condition, invertibility of Frechet derivative, in Theorem \ref{thm:3.1} is too artificial.
Now we turn out attention to get conditions to guarantee the invertibility of the Frechet derivative.
\begin{thm}\label{thm:3.2}
Suppose $(u,v)$ is a positive solution to (\ref{eq:ja}).\\
If $4\inf\frac{\partial g(x,y)}{\partial x}\inf\frac{\partial h(x,y)}{\partial y}uv > [\sup\frac{\partial g(x,y)}{\partial y}u + \sup\frac{\partial h(x,y)}{\partial x}v]^{2}$, then the Frechet derivative of (\ref{eq:ja}) at $(u,v)$ is invertible.
\end{thm} 
\pf The Frechet derivative at $(u,v)$ is
$$A = \left( \begin{array}{ll}
-\Delta + g(u,v) + u\frac{\partial g(u,v)}{\partial u} - a & u\frac{\partial g(u,v)}{\partial v}\\
v\frac{\partial h(u,v)}{\partial u} & -\Delta + h(u,v) + v\frac{\partial h(u,v)}{\partial v} - d
\end{array} \right).$$
We need to show that $N(A) = \{0\}$ by Fredholm alternative. If
$$\left\{ \begin{array}{l}
-\Delta \varphi + (g(u,v) + u\frac{\partial g(u,v)}{\partial u} - a)\varphi + \frac{\partial g(u,v)}{\partial v}u\psi = 0,\\
-\Delta \psi + \frac{\partial h(u,v)}{\partial u}v\varphi + (h(u,v) + v\frac{\partial h(u,v)}{\partial v} - d)\psi = 0, 
\end{array} \right.$$
then 
$$\left. \begin{array}{l}
\int_{\Omega}[|\nabla \varphi|^{2} + (g(u,v) + u\frac{\partial g(u,v)}{\partial u} - 
a)\varphi^{2} + \frac{\partial g(u,v)}{\partial v}u \varphi \psi ] = 0,\\
\int_{\Omega}[|\nabla \psi|^{2} + \frac{\partial h(u,v)}{\partial u}v\varphi\psi + (h(u,v) + v\frac{\partial h(u,v)}{\partial v} - d)\psi^{2}] = 0.
\end{array} \right.$$
Since $\lambda_{1}(g(u,v) - a) = \lambda_{1}(h(u,v) - d) = 0$,
$$\left. \begin{array}{l}
\int_{\Omega}[|\nabla \varphi|^{2} + (g(u,v) - a)\varphi^{2}] \geq 0,\\
\int_{\Omega}[|\nabla \psi|^{2} + (h(u,v) -d)\psi^{2}] \geq 0. 
\end{array} \right.$$     
Hence,
$$\left. \begin{array}{l}
\int_{\Omega}(u\frac{\partial g(u,v)}{\partial u}\varphi^{2} + \frac{\partial g(u,v)}{\partial v}u\varphi\psi) \leq 0,\\
\int_{\Omega}(\frac{\partial h(u,v)}{\partial u}v\varphi\psi + \frac{\partial h(u,v)}{\partial v}v\psi^{2}) \leq 0.
\end{array} \right.$$
Hence, $\int_{\Omega}[u\frac{\partial g(u,v)}{\partial u}\varphi^{2} + (\frac{\partial g(u,v)}{\partial v}u + \frac{\partial h(u,v)}{\partial u}v)\varphi\psi + \frac{\partial 
h(u,v)}{\partial v}v\psi^{2}] \leq 0.$\\
Hence, if $4\inf(\frac{\partial g(x,y)}{\partial x})\inf(\frac{\partial h(x,y)}{\partial y})uv > [(\sup(\frac{\partial g(x,y)}{\partial y}))u + (\sup(\frac{\partial h(x,y)}{\partial x}))v]^{2}$, then the integrand in the left side is positive definite form in $\Omega$, which means $\varphi \equiv \psi \equiv 0.$ Therefore, the above Frechet derivative $A$ is invertible. 
\\[0.2in]
Combining the Theorems \ref{thm:2.1}, \ref{thm:3.1} and \ref{thm:3.2}, we have the following which is actually the main result in this section.
\begin{cor}\label{cor:3.3} Suppose\\
$(A)$ $a > \lambda_{1} + g(0,k_{2})$, $d > \lambda_{1} + h(k_{1},0),$ and\\
$(B)$ 
$$\left.\begin{array}{lll}
4\inf_{B}\frac{\partial g(x,y)}{\partial x}\inf_{B}\frac{\partial h(x,y)}{\partial y} & > & [\sup\frac{\partial g(x,y)}{\partial y} + \sup\frac{\partial h(x,y)}{\partial x} \sup\frac{\theta_{d-h(0,\cdot)}}{\theta_{a-g(\cdot,k_{2})}}]\\
& & [\sup\frac{\partial g(x,y)}{\partial y} 
\sup\frac{\theta_{a-g(\cdot,0)}}{\theta_{d-h(k_{1},\cdot)}} + \sup\frac{\partial h(x,y)}{\partial x}],
\end{array}\right.$$
where $B = [0,k_{1}] \times [0,k_{2}]$.\\
Then there is a neighborhood $V$ of $(a,d)$ in $R^{2}$ such that if $(a_{0},d_{0}) \in V$, then 
(\ref{eq:ja}) with $(a,d) = (a_{0},d_{0})$ has a unique coexistence state.  
\end{cor}
\pf From $\theta_{a-g( ,0)} < k_{1},\;\theta_{d-h(0, )} < k_{2}$, and the monotonicity of $g(0,\cdot), h(\cdot,0)$ we have
$$\left\{ \begin{array}{l}
a > \lambda_{1} + g(0,k_{2}) \geq \lambda_{1}(g(0,\theta_{d-h(0,\cdot)})),\\ 
d > \lambda_{1} + h(k_{1},0) \geq \lambda_{1}(h(\theta_{a-g(\cdot,0)},0)).
\end{array} \right.$$
\begin{eqnarray*}
& & 4\inf_{B}\frac{\partial g(x,y)}{\partial x}\inf_{B}\frac{\partial h(x,y)}{\partial y}\\ 
& > & [ \sup\frac{\partial g(x,y)}{\partial y} + \sup\frac{\partial h(x,y)}{\partial x} \sup \frac{\theta_{d-h(0,\cdot)}}{\theta_{a-g(0,k_{2})}}] [\sup\frac{\partial g(x,y)}{\partial y} 
\sup\frac{\theta_{a-g(\cdot,0)}}{\theta_{d-h(k_{1},\cdot)}}\\
& & + \sup\frac{\partial h(x,y)}{\partial x}]\\
& = & [\sup\frac{\partial g(x,y)}{\partial y}]^{2}\sup\frac{\theta_{a-g(\cdot,0)}}{\theta_{d-h(k_{1},\cdot)}} + \sup\frac{\partial g(x,y)}{\partial y}\sup\frac{\partial h(x,y)}{\partial x}\\
& & + \sup\frac{\theta_{a-g(\cdot,0)}}{\theta_{d-h(k_{1},\cdot)}}\sup 
\frac{\theta_{d-h(0,\cdot)}}{\theta_{a-g(\cdot,k_{2})}}\sup\frac{\partial g(x,y)}{\partial y}\sup\frac{\partial h(x,y)}{\partial x}\\ 
& & + [\sup\frac{\partial h(x,y)}{\partial x}]^{2}\sup\frac{\theta_{d-h(0,\cdot)}}{\theta_{a-g(\cdot,k_{2})}}\\
& \geq & [\sup\frac{\partial g(x,y)}{\partial y}]^{2}\sup\frac{\theta_{a-g(\cdot,0)}}{\theta_{d-h(k_{1},\cdot)}} + 2\sup\frac{\partial g(x,y)}{\partial y}\sup\frac{\partial h(x,y)}{\partial x}\\
& & + [\sup\frac{\partial h(x,y)}{\partial x}]^{2}\sup\frac{\theta_{d-h(0,\cdot)}}{\theta_{a-g(\cdot,k_{2})}}
\end{eqnarray*}
since $\theta_{a-g(\cdot,0)} > \theta_{a-g(\cdot,k_{2})}, \theta_{d-h(0,\cdot)} > \theta_{d-h(k_{1},\cdot)}$.\\
Therefore, (\ref{eq:ja}) has a unique coexistence state $(u,v)$ from Theorem \ref{thm:2.1}. Furthermore, by the estimate of the solution in the proof of Theorem \ref{thm:2.1},
\begin{eqnarray*}
& & 4\inf_{B}\frac{\partial g(x,y)}{\partial x}\inf_{B}\frac{\partial h(x,y)}{\partial y}\\
& > & [\sup\frac{\partial g(x,y)}{\partial y} + \sup\frac{\partial h(x,y)}{\partial x} \sup\frac{\theta_{d-h(0,\cdot)}}{\theta_{a-g(\cdot,k_{2})}}] [\sup\frac{\partial g(x,y)}{\partial y}\sup\frac{\theta_{a-g(\cdot,0 )}}{\theta_{d-h(k_{1},\cdot)}}\\
& & + \sup\frac{\partial h(x,y)}{\partial x}]\\ 
& \geq & [\sup\frac{\partial g(x,y)}{\partial y} + \sup\frac{\partial h(x,y)}{\partial x} \frac{v}{u}][\sup\frac{\partial g(x,y)}{\partial y}\frac{u}{v} + \sup\frac{\partial h(x,y)}{\partial x}].
\end{eqnarray*} 
Thus, we obtain  
$$4\inf_{B}\frac{\partial g(x,y)}{\partial x}\inf_{B}\frac{\partial h(x,y)}{\partial y}uv > [\sup\frac{\partial g(x,y)}{\partial y}u + \sup\frac{\partial h(x,y)}{\partial x}v]^{2}.$$ 
It implies that the Frechet derivative of (\ref{eq:ja}) at $(u,v)$ is invertible from Theorem \ref{thm:3.2}. Therefore, the theorem follows from Theorem \ref{thm:3.1}.

\section{Uniqueness in a region of reproduction rates} 
Consider the model
\begin{equation}
\left\{ \begin{array}{l} 
\left. \begin{array}{l}
\Delta u + u(a - g(u,v)) = 0\\
\Delta v + v(d - h(u,v)) = 0
\end{array} \right. \mbox{in $\Omega$,}\\
u|_{\partial \Omega} = v|_{\partial \Omega} =0.
\end{array} \right. \label{eq:ga}
\end{equation}
Here $\Omega$ is a bounded smooth domain in $R^{n}$ and $g, h \in C^{1}$ are strictly increasing functions with respect to $u$ and $v$, and $g(0,0) = h(0,0) = 0$.
\\[0.1in]
The following is the main theorem.
\begin{thm}\label{thm:4.1} Suppose\\
$(A)$ $\Gamma$is a closed, convex region in $R^{2}$ such that for all $(a,d) \in \Gamma$,\\
\makebox[1.7in][r]{$a > \lambda_{1}(g(0,\theta_{d-h(0,\cdot)}))$} and $d > \lambda_{1}(h(\theta_{a-g(\cdot,0)},0))$,\\
$(B)$ there exist $c_{0} > 0$ and $c_{1} > 0$ such that for all $(a,d) \in \Gamma$, $g(x,0) > a > \lambda_{1}$ for $x > c_{0}$ and $h(0,y) > d > \lambda_{1}$ for $y > c_{1}$,\\
$(C)$ (\ref{eq:ga}) has a unique positive solution for every $(a,d) \in \partial_{L}\Gamma$,\\
\makebox[0.5in][r]{where} $\partial_{L}\Gamma = \{(\lambda_{d},d) \in \Gamma | \mbox{For any fixed}\;\;d, \lambda_{d} = \inf\{a | (a,d) \in \Gamma\}\}$,\\ 
$(D)$ for all $(a,d) \in \Gamma$, the Frechet derivative of (\ref{eq:ga}) at every positive 
solution\\
\makebox[0.4in][r]{to} (\ref{eq:ga}) is invertible.\\
Then for all $(a,d) \in \Gamma$, (\ref{eq:ga}) has a unique positive solution. Furthermore, there is an open set $W$ in $R^{2}$ such that $\Gamma \subseteq W$ and for every $(a,d) \in W$, (\ref{eq:ga}) has a unique positive solution.
\end{thm}
Theorem \ref{thm:4.1} goes even further than Theorem \ref{thm:3.1} which states the uniqueness in the whole region of $(a,d)$ whenever we have the uniqueness on the left boundary and invertibility of linearized operator at any particular solution inside the domain.
\pf For each fixed d, let $\lambda^{d} = \sup\{a : (a,d) \in \Gamma\}$ and $\lambda_{d} = \inf\{a | (a,d) \in \Gamma\}$. We need to show that for every $a$ such that $\lambda_{d} \leq a \leq \lambda^{d}$, (\ref{eq:ga}) has a unique positive solution. Since (\ref{eq:ga}) with $(a,d) = (\lambda_{d},d)$ has a unique positive solution $(u,v)$ and the Frechet derivative of (\ref{eq:ga}) at $(u,v)$ is invertible, by theorem \ref{thm:3.1}, there is an open neighborhood $V$ of $(\lambda_{d},d)$ in $R^{2}$ such that if $(a_{0},d_{0}) \in V$, then (\ref{eq:ga}) with $(a,d) = (a_{0},d_{0})$ has a unique positive solution. Let $\lambda_{s} = \sup\{\lambda \geq \lambda_{d} : \mbox{(\ref{eq:ga}) has a unique coexistence state for}\;\; \lambda_{d} \leq a \leq \lambda\}.$ We need to show that $\lambda_{s} \geq \lambda^{d}$. Suppose $\lambda_{s} < \lambda^{d}$. From the definition of $\lambda_{s}$, there is a sequence $\{\lambda_{n}\}$ such that $\lambda_{n} \rightarrow \lambda_{s}^{-}$ and there is a sequence $(u_{n},v_{n})$ of the unique positive solution of (\ref{eq:ga}) with $(a,d) =(\lambda_{n},d)$. Then by the Elliptic theory, there is $(u_{0},v_{0})$ such that $(u_{n},v_{n})$ converges to $(u_{0},v_{0})$ uniformly and $(u_{0},v_{0})$ is the solution to (\ref{eq:ga}) with $(a,d) = (\lambda_{s},d)$. We claim that $u_{0}$ is not identically zero and $v_{0}$ is not identically zero. Suppose this is false. Then by the Maximum Principle, one of the following cases should hold: $(1)$ $u_{0}$ is identically zero and $v_{0}$ is not identically zero. $(2)$ $u_{0}$ is not identically zero and $v_{0}$ is identically zero. $(3)$ Both $u_{0}$ and $v_{0}$ are identically zero. The argument is similar to what we had in the previous section.\\
$(1)$ Suppose $u_{0}$ is identically zero.\\
Let $\tilde{u_{n}} = \frac{u_{n}}{\parallel u_{n} \parallel_{\infty}}$ and $\tilde{v_{n}} = v_{n}$ for all $n \in N$. Then
$$\left\{ \begin{array}{l}
\Delta \tilde{u_{n}} + \tilde{u_{n}}(\lambda_{n} - g(u_{n},\tilde{v_{n}})) = 0\\
\Delta \tilde{v_{n}} + \tilde{v_{n}}(d - h(u_{n},\tilde{v_{n}})) = 0
\end{array} \right. \mbox{in $\Omega$.}$$
We know $\tilde{u_{n}} \rightarrow \tilde{u}$ from the elliptic theory, and
$$\left\{ \begin{array}{l}
\Delta \tilde{u} + \tilde{u} (\lambda_{s} - g(0,v_{0})) = 0\\
\Delta v_{0} + v_{0} (d - h(0,v_{0})) = 0
\end{array} \right. \mbox{in $\Omega,$}$$
since $g, h$ are continuous. Hence, $v_{0} = \theta_{d-h(0,\cdot)}$ and $\lambda_{s} = \lambda_{1}(g(0,v_{0}))$. If $v_{0}$ is identically zero, then by the monotonicity of $g$ and $\lambda_{1}$, we have $\lambda_{s} = \lambda_{1}(g(0,v_{0})) = \lambda_{1}(g(0,0)) \leq \lambda_{1}(g(0,\theta_{d-h(0,\cdot)})) < \lambda_{d}$ which is impossible. If $v_{0}$ is not identically zero, then $v_{0} = \theta_{d-h(0,\cdot)}$ and so $\lambda_{s} = \lambda_{1}(g(0,v_{0})) = \lambda_{1}(g(0,\theta_{d-h(0,\cdot)})) < \lambda_{d}$ which is also impossible.\\
$(2)$ Suppose $v_{0}$ is identically zero.\\
Let $\tilde{u_{n}} = u_{n}$ and $\tilde{v_{n}} = \frac{v_{n}}{\parallel v_{n} \parallel_{\infty}}$ for all $n \in N$. Then
$$\left\{ \begin{array}{l}
\Delta \tilde{u_{n}} + \tilde{u_{n}} (\lambda_{n} - g(\tilde{u_{n}},v_{n})) = 0\\
\Delta \tilde{v_{n}} + \tilde{v_{n}} (d - h(\tilde{u_{n}},v_{n})) = 0
\end{array} \right. \mbox{in $\Omega$.}$$
Again $\tilde{v_{n}} \rightarrow \tilde{v}$ by the elliptic theory, and
$$\left\{ \begin{array}{l}
\Delta u_{0} + u_{0} (\lambda_{s} - g(u_{0},0)) = 0\\
\Delta \tilde{v} + \tilde{v} (d - h(u_{0},0)) = 0
\end{array} \right. \mbox{in $\Omega$}$$
since $g, h$ are continuous. Hence, $d = \lambda_{1}(h(u_{0},0))$. If $u_{0}$ is identically zero, then by the monotonicity of $h$ and $\lambda_{1}$, we have $d = \lambda_{1}(h(u_{0},0)) = \lambda_{1}(h(0,0)) \leq \lambda_{1}(h(\theta_{\lambda_{s}-g(\cdot,0)},0)) < \lambda_{1}(h(\theta_{\lambda^{d}-g(\cdot,0)},0))$ which is impossible, since $(\lambda^{d},d) \in \Gamma$. If $u_{0}$ is not identically zero,, then $u_{0} = \theta_{\lambda_{s}-g(\cdot,0)}$ and so, $d = \lambda_{1}(h(u_{0},0)) = \lambda_{1}(h(\theta_{\lambda_{s}-g(\cdot,0)},0)) < \lambda_{1}(h(\theta_{\lambda^{d}-g(\cdot,0)},0))$ which is also impossible, since $(\lambda^{d},d) \in \Gamma$.
Consequently, $u_{0} > 0, v_{0} > 0$ in $\Omega$, that is, $(u_{0},v_{0})$ is a coexistence of (\ref{eq:ga}) with $(a,d) = (\lambda_{s},d)$. Since $(\lambda_{s},d) \in \Gamma$, by the assumption, the Frechet derivative of (\ref{eq:ga}) with $(a,d) = (\lambda_{s},d)$ at $(u_{0},v_{0})$ is invertible. Hence, by the Implicit Function Theorem, there is an open neighborhood $U$ of $\lambda_{s}$ and an open neighborhood $V$ of $(u_{0},v_{0})$ such that if $a \in U$, then (\ref{eq:ga}) has a unique coexistence state in $V$. But, by the definition of $\lambda_{s}$, there is a sequence $\{\lambda_{n}'\} \subseteq U$ such that $\lambda_{n}' \rightarrow \lambda_{s}^{+}$ and there is a sequence $\{(u_{n}',v_{n}')\}$ of the coexistence state of (\ref{eq:ga}) with $(a,d) = (\lambda_{n}',d)$ such that $(u_{n}',v_{n}') \notin 
V$ for all $n \in N$. By the Elliptic Theory again, $u_{n}' \rightarrow u_{0}', v_{n}' \rightarrow v_{0}'$ and from the same argument above, $(u_{0}',v_{0}') \notin V$ is also a coexistence of (\ref{eq:ga}) with $(a,d) = (\lambda_{s},d)$. Since $(\lambda_{s},d) \in \Gamma$, by the assumption again, the Frechet derivative of (\ref{eq:ga}) at $(u_{0}',v_{0}')$ is invertible. Hence, by the Implicit Function Theorem again, there is an open neighborhood $U'$ of $\lambda_{s}$ and an open neighborhood $V'$ of $(u_{0}',v_{0}')$ such that if $a \in U'$, then (\ref{eq:ga}) has a unique coexistence state in $V'$. Consequently, there are points in the left side of $\lambda_{s}$ such that (\ref{eq:ga}) has two different coexistence states. That is a contradiction to the definition of $\lambda_{s}$. Hence, $\lambda_{s} \geq \lambda^{d}$ and the first part of the theorem is proved. Furthermore, by the assumption, for each $(a,d) \in \Gamma$, the Frechet derivative of (\ref{eq:ga}) at the unique solution $(u,v)$ is invertible. Hence, Theorem \ref{thm:3.1} concludes that there is an open neighborhood $V_{(a,d)}$ of $(a,d)$ in $R^{2}$ such that if $(a_{0},d_{0}) \in V_{(a,d)}$, then (\ref{eq:ga}) with reproduction rates $(a_{0},d_{0})$ has a unique coexistence state. Let $W = \bigcup_{(a,d)\in\Gamma}V_{(a,d)}$. Then $W$ is an open set in $R^{2}$ such that $\Gamma \subseteq W$ and for each $(a_{0},d_{0}) \in W$, (\ref{eq:ga}) has a unique coexistence state. 
\\[0.1in]
Apparently, Theorem \ref{thm:4.1} generalizes Theorem \ref{thm:3.1} and consequently, we have the following which is actually the main conclusion in this section.
\\[0.1in]
\begin{cor}\label{cor:4.2} Suppose\\
$(A)$ $\Gamma$ is a closed, convex region in $R^{2}$,\\
$(B)$ there exist $k_{1}, k_{2} > 0$ such that for all $(a,d) \in \Gamma$, $a > \lambda_{1} + g(0,k_{2})$, $d > \lambda_{1} + h(k_{1},0), a - g(k_{1},0) < 0, d - h(0,k_{2}) < 0,$\\
$(C)$ 
$$\left.\begin{array}{lll}
4\inf_{B}\frac{\partial g(x,y)}{\partial x}\inf_{B}\frac{\partial h(x,y)}{\partial 
y} & > & [\sup\frac{\partial g(x,y)}{\partial y} + \sup\frac{\partial h(x,y)}{\partial x} \sup_{(a,d)\in\Gamma}\frac{\theta_{d-h(0,\cdot)}}{\theta_{a-g(\cdot,k_{2})}}]\\
& & [\sup\frac{\partial g(x,y)}{\partial y}\sup_{(a,d)\in\Gamma}\frac{\theta_{a-g(\cdot,0 )}}{\theta_{d-h(k_{1},\cdot)}} + \sup\frac{\partial h(x,y)}{\partial x}],
\end{array}\right.$$
where $B = [0,k_{1}] \times [0,k_{2}]$.\\
Then there is an open set $W$ in $R^{2}$ such that $\Gamma \subseteq W$ and for every $(a,d) \in \Gamma$, (\ref{eq:ga}) has a unique positive solution.   
\end{cor}
The condition $(B)$ means $\Gamma$ is some set of large self-reproduction rates, and the condition $(C)$ implies that the self-limitation rates are relatively larger than competition rates. Then the conclusion says that the existence and uniqueness of coexistence state are guaranteed on $\Gamma$ and the region $\Gamma$ can be extended to a larger set without losing the uniqueness.
\pf From $\theta_{a-g(\cdot,0)} < k_{1},\;\theta_{d-h(0,\cdot)} < k_{2}$, and the monotonicity of $g(0,\cdot), h(\cdot,0)$ we have
$$\left\{ \begin{array}{l}
a > \lambda_{1} + g(0,k_{2}) \geq \lambda_{1}(g(0,\theta_{d-h(0,\cdot)})),\\ 
d > \lambda_{1} + h(k_{1},0) \geq \lambda_{1}(h(\theta_{a-g(\cdot,0)},0)).
\end{array} \right.$$
for all $(a,d) \in \Gamma$.\\
By the condition $(C)$, for every $(a,d) \in \partial\Gamma$, 
\begin{eqnarray*}
& & 4\inf_{B}\frac{\partial g(x,y)}{\partial x}\inf_{B}\frac{\partial h(x,y)}{\partial y}\\ & > & [ \sup\frac{\partial g(x,y)}{\partial y} + 
\sup\frac{\partial h(x,y)}{\partial x} \sup \frac{\theta_{d-h(0,\cdot)}}{\theta_{a-g(0,k_{2})}}] 
[\sup\frac{\partial g(x,y)}{\partial y}\sup\frac{\theta_{a-g(\cdot,0)}}{\theta_{d-h(k_{1},\cdot)}}\\
& & + \sup\frac{\partial h(x,y)}{\partial x}]\\
& = & [\sup\frac{\partial g(x,y)}{\partial y}]^{2}\sup\frac{\theta_{a-g(\cdot,0)}}{\theta_{d-h(k_{1},\cdot)}} + \sup\frac{\partial g(x,y)}{\partial y}\sup\frac{\partial h(x,y)}{\partial x}\\
& & + \sup\frac{\theta_{a-g(\cdot,0)}}{\theta_{d-h(k_{1},\cdot)}}\sup 
\frac{\theta_{d-h(0,\cdot)}}{\theta_{a-g(\cdot,k_{2})}}\sup\frac{\partial g(x,y)}{\partial y}\sup\frac{\partial h(x,y)}{\partial x}\\ 
& & + [\sup\frac{\partial h(x,y)}{\partial x}]^{2}\sup\frac{\theta_{d-h(0,\cdot)}}{\theta_{a-g(\cdot,k_{2})}}\\
& \geq & [\sup\frac{\partial g(x,y)}{\partial y}]^{2}\sup\frac{\theta_{a-g(\cdot,0)}}{\theta_{d-h(k_{1},\cdot)}} + 2\sup\frac{\partial g(x,y)}{\partial y}\sup\frac{\partial h(x,y)}{\partial x}\\
& & + [\sup\frac{\partial h(x,y)}{\partial x}]^{2}\sup\frac{\theta_{d-h(0,\cdot)}}{\theta_{a-g(\cdot,k_{2})}}
\end{eqnarray*}
since $\theta_{a-g(\cdot,0)} > \theta_{a-g(\cdot,k_{2})}, \theta_{d-h(0,\cdot)} > \theta_{d-h(k_{1},\cdot)}$.
Therefore, by the Theorem \ref{thm:2.1}, (\ref{eq:ga}) has a unique coexistence state for all $(a,d) \in \partial\Gamma$. Furthermore, by the estimate of the solution in the proof of Theorem \ref{thm:2.1}, if $(u,v)$ is a positive solution for $(a,d) \in\Gamma$, then
\begin{eqnarray*}
& & 4\inf_{B}\frac{\partial g(x,y)}{\partial x}\inf_{B}\frac{\partial h(x,y)}{\partial y}\\
& > & [\sup\frac{\partial g(x,y)}{\partial y} + \sup\frac{\partial h(x,y)}{\partial x} \sup\frac{\theta_{d-h(0,\cdot)}}{\theta_{a-g(\cdot,k_{2})}}] [\sup\frac{\partial g(x,y)}{\partial y}\sup\frac{\theta_{a-g(\cdot,0 )}}{\theta_{d-h(k_{1},\cdot)}}\\
& & + \sup\frac{\partial h(x,y)}{\partial x}]\\ 
& \geq & [\sup\frac{\partial g(x,y)}{\partial y} + \sup\frac{\partial h(x,y)}{\partial x} \frac{v}{u}][\sup\frac{\partial g(x,y)}{\partial y} \frac{u}{v} + \sup\frac{\partial h(x,y)}{\partial x}].
\end{eqnarray*} 
Thus, we obtain  
$$4\inf_{B}\frac{\partial g(x,y)}{\partial x}\inf_{B}\frac{\partial h(x,y)}{\partial y}uv > [\sup\frac{\partial g(x,y)}{\partial y}u + \sup\frac{\partial h(x,y)}{\partial x}v]^{2}.$$ 
It implies that if $(u,v)$ is a positive solution of (\ref{eq:ga}) for $(a,d) \in\Gamma$, then the Frechet derivative of (\ref{eq:ga}) at $(u,v)$ is invertible from Theorem \ref{thm:3.2}. 
Therefore, the theorem follows from Theorem \ref{thm:4.1}.
\\[0.2in]
ACKNOWLEDGMENT. The author would like to express his gratitude to his thesis advisor Professor Zheng Fang Zhou and his colleague Professor Shandelle Henson for their kind help and encouragement.

\end{document}